\documentclass[11pt, reqno]{amsart}
\usepackage{amsmath,amsthm,amssymb}
\theoremstyle{plain}
\newtheorem{thm}{Theorem}[section]
\newtheorem{prop}[thm]{Proposition}
\newtheorem{lem}[thm]{Lemma}
\newtheorem{cor}[thm]{Corollary}

\theoremstyle{definition}

\newtheorem{rem}[thm]{Remark}
\newtheorem{defn}[thm]{Definition}
\newtheorem{eg}[thm]{Example}
\newtheorem{subtitle}[thm]{}
\newtheorem{ex}{Exercise}[section]
\numberwithin{equation}{section}

\def\a{\alpha}

\def\d{\delta}
\def\D{\triangle}
\def\e{\epsilon}
\def\g{\gamma}
\def\G{\Gamma}
\def\K{\nabla}
\def\l{\lambda}
\def\L{\Lambda}
\def\n{\,\vert\,}
\def\o{\theta}

\def\W{\Omega}
\def\ca{{\mathcal{A}}}
\def\cd{{\mathcal{D}}}
\def\cf{{\mathcal{F}}}
\def\cg{{\mathcal{G}}}

\def\ck{{\mathcal{K}}}

\def\co{{\mathcal{O}}}
\def\cp{{\mathcal{P}}}
\def\cs{{\mathcal{S}}}

\def\cu{{\mathcal{U}}}

\def\li{\langle}
\def\ri{\rangle}
\def\n{\ \vert\ }

\def\bs{\bigskip}

\def\ss{\smallskip}

\def\ni{\noindent}
\def\ti{\tilde}
\def\p{\partial}

\def\Im{{\rm Im\/}}
\def\I{{\rm I\/}}
\def\II{{\rm II\/}}
\def\diag{{\rm diag}}

\def\Gr{{\rm Gr}}

\def\R{\mathbb{R} }
\def\C{\mathbb{C}}

\newcommand{\beg}{\begin{eg}}
\newcommand{\eeg}{\end{eg}}
\newcommand{\bthm}{\begin{thm}}
\newcommand{\ethm}{\end{thm}}
\newcommand{\bprop}{\begin{prop}}
\newcommand{\eprop}{\end{prop}}
\newcommand{\bcor}{\begin{cor}}
\newcommand{\ecor}{\end{cor}}
\newcommand{\blem}{\begin{lem}}
\newcommand{\elem}{\end{lem}}
\newcommand{\bca}{\begin{cases}}
\newcommand{\eca}{\end{cases}}
\newcommand{\brem}{\begin{rem}}
\newcommand{\erem}{\end{rem}}
\newcommand{\bpm}{\begin{pmatrix}}
\newcommand{\epm}{\end{pmatrix}}
\newcommand{\bbm}{\begin{bmatrix}}
\newcommand{\ebm}{\end{bmatrix}}
\newcommand{\bvm}{\begin{vmatrix}}
\newcommand{\evm}{\end{vmatrix}}
\newcommand{\bdefn}{\begin{defn}}
\newcommand{\edefn}{\end{defn}}
\newcommand{\bsub}{\begin{subtitle}}
\newcommand{\esub}{\end{subtitle}}
\newcommand{\bex}{\begin{ex}}
\newcommand{\eex}{\end{ex}}
\newcommand{\ben}{\begin{enumerate}}
\newcommand{\een}{\end{enumerate}}

\def\half{\frac{1}{2}}
\def\pt{\frac{\p}{\p t}}
\def\px{\frac{\p}{\p x}}
\def\py{\frac{\p}{\p y}}

\def\pzbar{\frac{\p}{\p \bar z}}

\def\unon{\frac{U(n)}{O(n)}}

\begin{document}

\title[Geometric soliton equations]
{Applications of Loop Group Factorization to Geometric Soliton Equations}
\author{{Chuu-Lian Terng}$^*$}\thanks{$^*$Research supported in   part by NSF grant DMS- 052975}
\address{Department of Mathematics \\
University of California at Irvine, Irvine, CA 92697-3875}
\email{cterng@math.uci.edu}

\begin{abstract}

The $1$-d Schr\"odinger flow on $S^2$, the Gauss-Codazzi equation for flat Lagrangian submanifolds in $\R^{2n}$, and the space-time monopole equation are all examples of geometric soliton equations.  The linear systems with a spectral parameter (Lax pair) associated to these equations satisfy the reality condition associated to $SU(n)$.  In this article, we explain the method developed jointly with K. Uhlenbeck, that uses various loop group factorizations to construct inverse scattering transforms, B\"acklund transformations, and solutions to Cauchy problems for these equations.

\end{abstract}

\maketitle


\section{Introduction}

A Hamiltonian system in $2n$-dimension is called {\it completely integrable\/} if it has $n$ independent commuting Hamiltonians.  By the Arnold-Liouville Theorem, such systems have {\it action-angle\/} variables that linearize the flow.  The concept of completely integrability has been extended to  soliton equations. These equation can be  linearized using ``scattering data'', allowing one to use the Inverse Scattering method to solve the Cauchy problem with rapidly decaying initial data.   Two model examples are the Korteweg-de Vries equation (KdV) and the non-linear Sch\"rodinger equation (NLS). Soliton equations often arise naturally in differential geometry too. For example, the Gauss-Codazzi equations for surfaces in $\R^3$ with Gaussian curvature $-1$,  isothermic surfaces in $\R^3$ \cite{CieGolSym95}, isometric immersions of space forms in space forms \cite{FerPed96b, Ter97, Ten98}, Egoroff metrics, and  flat Lagrangian submanifolds in $\C^n$ and $\C P^n$ \cite{Ter03}, and the space-time monopole equation are soliton equations. 
  
  One of the key properties of a soliton equation is the existence of a Lax pair.  
A PDE for maps $q:\R^n\to \R^m$ admits a {\it Lax pair\/} if there exists a family of $\cg$-connections $\o_\l$ on $\R^n$, given in terms of $q$, such that the condition for $\o_\l$ to be flat for all $\l$ in an open subset of $\C$ is that $q$ satisfy the PDE. The parameter $\l$ is called the {\it spectral parameter\/}.  For a solution $q$ decaying at spatial infinity, we often can find a normalization so that there exists a unique parallel frame $E_\l$ of $\o_\l$.  Usually $E_\l$  has two types of singularities for $\l\in \C P^1$; one type is a jump across a contour and the other type is a pole.  We call the jump singularities of $E_\l$ the {\it continuous scattering data for $q$\/} and the poles and residues of $E_\l$ the {\it discrete scattering data for $q$\/}.  The {\it scattering transform\/} maps a solution $q$ to its scattering data $S$.  A key feature of soliton PDEs is that the induced equations on the scattering data is linear, so it is easy to write down the scattering data of a solution at time $t$ for a given initial data.  The inverse scattering transform  reconstructs $q$ from the scattering data,  i.e.,  it reconstructs $E_\l$ from prescribed singularities.  This is done for KdV in \cite{GGKM67}, for NLS in \cite{ZakMan73, FadTak87}, and for the $n$-wave equation in \cite{ZakMan73, BeaCoi84, BeaCoi85}.  As a consequence, the Cauchy problem for these soliton equations can be solved via the inverse scattering transform.  

The proof of the existence of the inverse scattering transform for soliton equations involves hard analysis and is difficult (\cite{BeaCoi84}).  However, if the Lax pair satisfies the $SU(n)$-reality condition, then the frame $E_\l(x)$ has only one jumping circle plus pole singularities in the $\l$- sphere for each $x$, so the continuous scattering data is a loop into $SL(n,\C)$ for each $x$.   In this case, we can use Pressley-Segal  loop group factorization  to construct the inverse scattering transform for the continuous scattering data (cf. \cite{TeUh98}). 

B\"acklund transformations (BTs) for surfaces in $\R^3$ with $K=-1$ arose from the study of  line congruences in classical differential geometry.  It associates to each surface in $\R^3$ with $K=-1$ a family of compatible systems of ordinary differential equations (ODEs) so that solutions of these ODE systems give rise to a family of new surfaces in $\R^3$ with $K=-1$. One can use line and sphere congruences to construct B\"acklund type transformations for many  geometric problems in  differential geometry (cf. \cite{Ten98}).  
B\"acklund transformations  for soliton equations produce a new solution from a given one by adding discrete scattering data.  These transformations can be obtained in a unified way from the following type of factorization:  Let $\G_1, \G_2$ be disjoint subsets of $S^2$, and $g_i:S^2\setminus \G_i\to GL(n,\C)$ holomorphic for $i=1, 2$.  Factor $g_1 g_2= \ti g_2 \ti g_1$ such that $\ti g_i$ is holomorphic on $S^2\setminus \G_i$.  This factorization can always be done when $g_1$ is rational and $g_i$ satisfy the $SU(n)$-reality condition, so global B\"acklund transformations exist for flows in the $SU(n)$-hierarchy and for the space-time monopole equation with gauge group $SU(n)$.  Moreover, if the initial data $q_0$ has continuous scattering data $S$ and discrete scattering data $\D$, then we can first use Pressely-Segal loop group factorization to construct a solution $q$ whose scattering data is $S$, and then apply BTs to $q$ to construct the solution $\ti q$ with scattering data $S\cup \D$.    
 
This paper is organized as follows:  In section 2, we outline the construction of the ZS-AKNS hierarchy of soliton equations associated to a complex simple Lie algebra $\cg$, and review certain invariant submanifolds and restricted flows associated to involutions of $\cg$.  We give examples of PDEs in submanifold geometry that are soliton equations in section 3. In section 4, we give a brief review of Lax pairs associated to the space-time monopole equations.  The direct scattering for soliton equations in the $SU(n)$-hierarchy is given in section 5,  and direct scattering  for  space monopole equation is in section 6.  We use Pressley-Segal loop group factorization to construct the inverse scattering transform for flows in the $SU(n)$-hierarchy and for the space-time monopole equation in section 7 and 8 respectively.  In section 9, we use the Birkhoff factorization to construct local solutions for flows in the $SU(n)$-hierarchy.  Finally, we discuss the constructions of B\"acklund transformations, pure solitons, and solutions with both continuous and discrete scattering data for flows in the $SU(n)$-hierarchy and for the space-time monopole equations in the last two sections. 

{\bf Acknowledgment:} The author  thanks her long-time collaborator and good friend Karen Uhlenbeck.  Much of this article concerns our joint project on the differential geometric aspects of soliton equations.

\bs
\section{Soliton equations associated to simple Lie algebras}

The method of constructing a hierarchy of $n\times n$ soliton flows developed 
by Zakharov-Shabat \cite{ZakSha79} and Ablowitz-Kaup-Newell-Segur \cite{AKNS74} works equally well if we replace the algebra of $n\times n$ matrices by a semi-simple, complex
Lie algebra $\cg$ (cf. \cite{HarSaiShn84b, Sat84, TeUh98}).   

\ss\ni {\bf The $G$-hierarchy}

\ss
Let $G$ be a complex, simple Lie group, $\cg$ its Lie algebra,  $\li\ , \ \ri$ a non-degenerate, ad-invariant bilinear form on $\cg$, $\ca$ a maximal abelian subalgebra of $\cg$,  and $\ca^\perp=\{\xi\in \cg\n \li \xi, \ca\ri=0\}$.  
Let $\cs(\R,\ca^\perp)$ denote the space  of rapidly decaying maps from $\R$ to $\ca^\perp$.
Fix a regular element   $a\in \ca$ (i.e., the centralizer $\cg_a=\ca$). Then
there is a unique family of $\cg$-valued maps $Q_{b,j}(u)$ parametrized by $b\in \ca$
and  positive integer $j$ satisfying the following recursive formula,
\begin{equation} \label{bi}
(Q_{b,j}(u))_x + [u, Q_{b, j}(u)] = [Q_{b,j+1}(u),a],
\quad Q_{b,0}(u)=b,
\end{equation}
 and $\sum_{j=0}^\infty Q_{b,j}(u)\l^{-j}$  is conjugate to $b$ as an  asymptotic expansion at $\l=\infty$.
In fact, $Q_{b,j}(u)$ is a polynomial in
$u, \p_x u, \ldots, \p_x^{j-1}u$ (cf., \cite{Sat84, TeUh98}).   For $b\in \ca$ and a positive integer $j$, the {\it $(b,j)$-flow\/}  is the following evolution equation on $\cs(\R, \ca^\perp)$:
\begin{equation} \label{co}
u_t= (Q_{b,j}(u))_x+ [u,Q_{b,j}(u)]=[Q_{b,j+1}(u),a].
\end{equation}
 The {\it $G$-hierarchy\/} is the collection of these $(b,j)$-flows.

 The recursive formula \eqref{bi} implies that  $u$ is a solution of the
$(b,j)$-flow \eqref{co} 
if and only if 
\begin{equation} \label{ex}
\o_\l=(a\l + u)\ dx + (b\l^j + Q_{b,1}(u)\l^{j-1} +
\cdots + Q_{b,j}(u))\ dt
\end{equation}
is a flat $\cg$-valued connection 1-form on the $(x,t)$ plane
for all $\l\in \C$. Here $\o_\l$ is (left) flat, i.e., $d\o_\l+ \o_\l\wedge \o_\l=0$.   In other words, $\o_\l$ is a Lax pair for the $(b,j)$-flow \eqref{co}.  Also $\o_\l$ is flat is equivalent to 
$$[\p_x + a\l + u, \ \p_t + b\l^j + Q_{b,1}(u) \l^{j-1} + \cdots + Q_{b,j}(u)]=0.$$

\ss
\ni {\bf The $U$-hierarchy}  \ss

Let $\tau$ be an involution of $G$ such that its differential at the identity $e$ (still denoted by
$\tau$) is a conjugate linear involution on the complex Lie algebra $\cg$, and $U$ the fixed point set of
$\tau$.  The Lie algebra $\cu$ of $U$ is a real form of
$\cg$.  If $a, b\in \cu$, then the
$(b,j)$-flow in the $G$-hierarchy leaves $\cs(\R, \ca^\perp\cap \cu)$-invariant (cf.  
\cite{TeUh98}). The restriction of the flow \eqref{co} to $\cs(\R, \ca^\perp\cap \cu)$ is the
$(b,j)$-flow in the {\it $U$-hierarchy\/}. Since $Q_{b,j}(u)$ lies in $\cu$,  the Lax pair $\o_\l$ defined by \eqref{ex} is a $\cg$-valued
$1$-form satisfying the {\it
$U$-reality condition\/}:
\begin{equation}\label{fp}
\tau(\o_{\bar\l})= \o_\l.
\end{equation}

\ni{\bf The $U/K$-hierarchy\/}
\ss

Suppose $U$ is the real form defined by the involution $\tau$ of $G$, and $\sigma$ an involution of $G$ such that $d\sigma_e$ is complex linear and $\sigma\tau=\tau\sigma$.  Let $K$ be the fixed point set of $\sigma$ in $U$, $\cu$ and $\ck$ the Lie algebras of $U$ and $K$ respectively, and $\cp$ the $-1$-eigenspace of $d\sigma_e$ on $\cu$.  Then $U/K$ is a symmetric space, and $\cu=\ck+\cp$.   Let $\ca$ be a maximal abelian subalgebra in $\cp$.  If $a, b\in \ca$ and $u$ in $\ca^\perp\cap \ck$, then the $(b,j)$-flow in the $U$-hierarchy leaves $\cs(\R, \ca^\perp\cap \ck)$ invariant if $j$ is odd, and is normal to $\cs(\R, \ca^\perp\cap \ck)$ if $j$ is even.  The restriction of odd flows in the $U$-hierarchy  to $\cs(\R, \ca^\perp\cap \ck)$ is called the {\it $U/K$-hierarchy\/}.  Moreover, $\o_\l$ satisfies the {\it $U/K$-reality condition\/}
$$\tau(\o_{\bar \l})= \o_\l, \quad \sigma(\o_\l)=\o_{-\l}.$$

\begin{eg}  {\it $SL(2,\C)$-hierarchy\/} (cf. \cite{AblCla91}).  

 Let  $a=b =\diag(i,-i)$, and  $\ca=\C a$. Then
$$\ca^\perp=\left\{\begin{pmatrix}
0&q\cr r&0\cr
\end{pmatrix} 
\bigg| \  q,r\in \C\right\},$$
\begin{align*}
& Q_{a,1}(u)= u= \begin{pmatrix}
0&q\cr r&0 \end{pmatrix}, \quad  Q_{a,2}(u)=\frac{i}{2}\begin{pmatrix}
qr& \p_xq\cr -\p_x r&
-qr\end{pmatrix}, \\
& Q_{a,3}=\frac{i}{4}\bpm q\p_x r - r\p_x q& -\p_x^2 q + 2 q^2r\\ -\p_x^2 r+ 2 qr^2& -q\p_x r + r\p_x q\epm, \quad  \cdots.
\end{align*}
The $(a,j)$-flows, $j=1,2, 3$, in the $SL(2,\C)$-hierarchy are:
\begin{align*}
&\p_t q=\p_xq, \ \ \p_t r= \p_x r,\ \ \cr
&\p_tq= \frac{i}{2} (\p_x^2q - 2q^2 r), \ \ \p_t r= -\frac{i}{2} (\p_x^2 r - 2qr^2),\\
& \p_t q= \frac{1}{4}(-\p_x^3 q+ 6qr\p_x q),\ \ \p_t r= \frac{1}{4}(-\p_x^3 r + 6 qr\p_x r).
\end{align*}

Let $\tau$ be the involution of $sl(2,\C)$ defined by
$\tau(\xi)= -\bar\xi^t$. Then the fixed point set of $\tau$ is the real form $\cu= su(2)$ and $$\ca^\perp\cap \cu=\left\{\begin{pmatrix}
0&q\cr -\bar q&0\cr
\end{pmatrix}
\bigg| \ q\in \R\right\}.$$ 
So the $SU(2)$-hierarchy is the
restriction of the $SL(2,\C)$-hierarchy to the subspace $r=-\bar q$.   The second flow in the
$SU(2)$-hierarchy is the NLS
$\p_t q= \frac{i}{2} (\p_x^2 q + 2 |q |^2 q) $.

Let $\sigma(\xi)= -(\xi^t)$. Then $\sigma\tau=\tau\sigma$ and the corresponding symmetric space is $SU(2)/SO(2)$.  Note that $u\in \cs(\R, \ca^\perp\cap \ck)$ means $q=-r$ is real.  The third flow in the $SU(2)/SO(2)$-hierarchy is the mKdV equation $q_t=- \frac{1}{4}(q_{xxx} + 6 q^2 q_x)$.  

\eeg

\ss
\ni{\bf The $U$-system}
\ss

Let $\cu$ be the real form of $\cg$ defined by the involution $\tau$,  $\ca$ a maximal abelian subalgebra of $\cu$, and $a_1, \ldots, a_n$ a basis of $\ca$. The $U$-system is the following PDE for $v: \R^n\to \ca^\perp$: 
 \begin{equation}\label{ag}
 [a_j, \p_{x_i} v] - [a_i, \p_{x_j} v] + [[a_i, v], [a_j, v]]=0, \quad i\not=j.
 \end{equation}
 It has a Lax pair 
 \begin{equation}\label{cf}
 \o_\l = \sum_{i=1}^n(a_i\l + [a_i, v])\ dx_i.
\end{equation}
This  Lax pair satisfies the {\it $U$-reality condition\/} $\o_\l= \tau(\o_{\bar\l})$. 
 
 \ss\ni{\bf The $U/K$-system}
 \ss
 
 Let $\tau, \sigma, U, K, \cp, \ca$ be as in the $U/K$-hierarchy, and  $a_1, \ldots, a_n$ a basis of $\ca$.  The $U/K$-system is the restriction of \eqref{ag} to the space of $v:\R^n\to \ca^\perp\cap \cp$. Since $a_i\in \cp$ and $[a_i, v]\in \ck$, its Lax pair $\o_\l = \sum_{i=1}^n (a_i \l + [a_i, v])\ dx_i$ satisfies the {\it $U/K$-reality condition\/}.
 
 \ss
\ni{\bf The frame of a Lax pair}
\hfil
\ss

Given a family of flat $\cg$-valued connections $\o_\l=\sum_{i=1}^n P_i(x,\l) dx_i$ on $\R^n$, we call $E(x, \l)$ {\it a frame\/}  of $\o_\l$ if
$E^{-1}\p_{x_i} E= P_i$ for all $ 1\leq i\leq n$. 

\bprop\label{ci}
Let $G, \tau, \sigma, U$ and $K$ be above, and $E_\l$ the frame of $\o_\l$ such that $E_\l(0)=\I$.  
\ben
\item If $\o_\l$ satisfies the $U$-reality condition, then $E_\l$ satisfies the $U$-reality condition $\tau(E_{\bar\l})=E_\l$,
\item If $\o_\l$ satisfies the $U/K$-reality condition, then $E_\l$ satisfies the $U/K$-reality condition
$\tau(E_{\bar \l})=E_\l, \ \  \sigma(E_\l)= E_{-\l}$.
\een
\eprop

\bs
\section{Soliton equations in submanifold geometry}

Since the Gauss-Codazzi equations for submanifolds in space forms are equivalent to the the flatness of certain connections, it is not surprising that many PDEs in submanifold geometry turns out to be soliton equations.  We give some examples below:

\beg[\bf Vortex filament equation, Schr\"odinger flow on $S^2$, and the NLS]\hfil

In 1906, da Rios modeled the movement of a thin vortex in a viscous fluid
by the motion of a curve propagating in $\R^3$ by
\begin{equation} \label{aa}
\p_t\g= \p_x \g\times \p_x^2\g.
\end{equation}  
If $\g$ is a solution of \eqref{aa}, then 
$$\p_t\li \p_x \g, \p_x\g\ri = 2\li \p_x\p_t \g, \p_x\g\ri = 2\li \p_x(\p_x\g \times\p_x^2 \g, \p_x\g\ri =0.$$  So \eqref{aa} preserves arc-length.  Hence we may assume that  a solution $\g(x,t)$ of \eqref{aa} satisfying $|| \p_x \g||=1$.  It is known that 
there exists a parallel normal frame $(v_1, v_2)(\cdot, t)$ for each curve $\g(\cdot, t)$ such that $q= k_1+ ik_2$ is a solution of the NLS, where $k_1$ and $k_2$  are the principal curvatures of $\g$ along $v_1$ and $v_2$ respectively.

Let $\mathcal{E}$ denote the energy functional on the space of paths on $S^2$, and $J$ the complex structure on $S^2$ (if we view $S^2\subset \R^3$, then $J_u(v)= u\times v$).  The Schr\"odinger flow on $S^2$ is 
$$u_t= J_u(\K {\mathcal{E}}(u)) = u\times u_{xx}.$$
 If $\g$ is a solution of \eqref{aa}, then $u= \g_x$ is a solution of the Schr\"odinger flow on $S^2$ (\cite{FadTak87,TeUh99}).  

\eeg

\beg[\bf Isothermic surfaces in $\R^3$]\hfil

A  parametrized surface $f(x,y)\in \R^3$ is {\it isothermal\/} if $(x,y)$ is a conformal line of curvature coordinate system, i.e., the two fundamental forms are of the form
$$\I= e^{2u} (dx_1^2+ dx_2^2), \quad \II= e^{u}(r_1 dx_1^2 + r_2 dx_2^2).$$
The Gauss-Codazzi equation is the $\frac{O(4,1)}{O(3)\times O(1,1)}$-system (cf. \cite{CieGolSym95,  Bur00, BrDuPaTe02}).   
\eeg

\beg[\bf Local isometric immersions of $N^n(c)$ in $N^{2n}(c)$]\label{cb}\cite{Ter97}\hfil

Let $N^m(c)$ denote the $n$-dimensional space form of constant sectional curvature $c$.  The normal bundle of a submanifold $M$ in $N^m(c)$ is {\it flat\/} if its induced normal connection is flat, and is {\it non-degenerate\/} if the dimension of $\{A_v\n v\in \nu(M)_p\}$ is equal to ${\rm codim\/}(M)$. Here $A_v$ is the shape operator along normal vector $v$.   It is proved in \cite{Ter97} that if $M^n$ is a submanifold of $N^{2n}(c)$ with constant sectional curvature $c$ and its normal bundle $\nu(M)$ is flat and non-degenerate, then there exists a local orthogonal coordinate system $(x_1, \ldots, x_n)$ on $M$ and parallel normal frame $e_{n+1}, \ldots, e_{2n}$ such that
\begin{equation}\label{ae}
\I=\sum_{i=1}^n b_i^2 dx_i^2,\quad \II= \sum_{j=1}^n a_{ji}b_idx_i^2 e_{n+j}.
\end{equation}
Moreover, the Levi-Civita connection $1$-form for $\I$ is $w= \d F - F^t\d$, where $F=(f_{ij})$, $f_{ij}= \frac{\p_{x_j} b_i}{b_j}$ if $ i\not=j$, $f_{ii}=0$ for all $1\leq i\leq n$,   and $\d= \diag(dx_1, \ldots, dx_n)$.
The Gauss-Codazzi equation for the local isometric immersion becomes an equation for $F$, which is  the $\frac{O(2n)}{O(n)\times O(n)}$-system if $c=0$, the $\frac{O(2n+1)}{O(n+1)\times O(n)}$-system if $c=1$, and the $\frac{O(2n,1)}{O(n)\times O(n,1)}$-system if $c=-1$. 
\eeg

\beg[\bf  Egoroff metrics and the $\unon$-system]\hfil

 A local orthogonal system $(x_1, \ldots, x_n)$ of $\R^n$ is {\it Egoroff\/} if the flat Euclidean metric $ds^2$ written in this coordinate system is of the form 
$$ds^2= \sum_{j=1}^n \p_{x_i} \phi \ dx_i^2$$ for some smooth function $\phi$.  Then $F=(f_{ij})$ is a solution of the $\frac{U(n)}{O(n)}$-system, where
$f_{ij}=  \frac{\p_{x_i}\p_{x_j} \phi}{2\ \p_{x_i}\phi \p_{x_j} \phi}$ if $i\not= j$ and $f_{ii}=0$ for $1\leq i\leq n$.  
Conversely, given a solution $F=(f_{ij}):\R^n\to V_n$ of the $\unon$-system, the first order system 
\begin{equation}\label{cj}
\p_{x_j} b_i = f_{ij} b_j, \quad i\not=j
\end{equation}
is solvable for $b_1, \ldots, b_n$, and solutions are locally defined and depend on $n$ functions of one variables.  Moreover, since $f_{ij}= f_{ji}$, $\sum_{i=1}^n b_i^2 dx_i$ is closed, hence locally there exists a smooth function $\phi$ such that $b_i^2 = \p_{x_i}\phi$ for $1\leq i\leq n$.

Although we can construct global solutions $F$ for the $U(n)/O(n)$-system, it is not clear whether there exist global  solutions $b_i$ of \eqref{cj} such that $b_i>0$ and the metric $ds^2= \sum_{i=1}^n b_i^2$ is complete.  This is also the case for isometric immersions of $N^n(c)$ in $N^{2n}(c)$ and for the next example.

\eeg

\beg[\bf Flat Lagrangian submanifolds in $\R^{2n}$]\hfil

As seen in Example \ref{cb}, the Gauss-Codazzi equation for  local isometric immersions of $\R^n$ into $\R^{2n}$ with flat and non-degenerate normal bundle is the $\frac{O(2n)}{O(n)\times O(n)}$-system.  These immersions are Lagrangian if and only if $F$ is symmetric and $F$ is a solution of the $\unon$-system.  

\eeg

\bs
\section{The space-time monopole equation}

For flows in the $SU(n)$-hierarchy, we have been using left flat connections $\o=\sum_{i=1}^nA_idx_i$, i.e., $d\o+\o\wedge \o=0$ or equivalently, $[\p_{x_i}+A_i, \p_{x_j}+A_j]=0$ for all $i\not=j$.  But for space-time monopole equations,  it is more customary to use right flat connections, i.e., $d\o-\o\wedge \o=0$, or equivalently, $[\p_{x_i}-A_i, \p_{x_j}-A_j]=0$ for all $i\not=j$. 

The curvature of a $su(n)$-valued connection $1$-form $A=\sum_{i=1}^4 A_i(x) dx_i$ is 
$F_A=\sum_{i< j} F_{ij} \ dx_i\wedge dx_j$. where 
$$F_{ij}= [\p_{x_i} - A_i, \p_{x_j} -A_j] = \p_{x_j} A_i - \p_{x_i} A_j + [A_i, A_j].$$
The connection $A$ is anti self-dual Yang-Mills (ASDYM) if 
$$F_A= -\ast F_A,$$
where $\ast$ is the Hodge star operator with respect to the metric $dx_1^2+dx_2^2 -dx_3^2 - dx_4^2$.

Set $ z=x_1+ix_2, \quad w= x_3+ix_4$,
$\K_z= \frac{1}{2}(\K_1-i\K_2)= \frac{\p}{\p z} - A_z$, $\K_{\bar
z}= \frac{1}{2} (\K_1+i\K_2)=\pzbar -A_{\bar z}$, and $\K_w,
\K_{\bar w}$ similarly. Since $A_i\in u(n)$,  $A_{\bar z} =
-A_z^*$ and $A_{\bar w}=-A_w^*$.
Then (cf. \cite{BelZak78, MaZa81}) $A$ is ASDYM if and only if 
\begin{equation}\label{am}
[\K_{\bar w} +\mu \K_z, \ \K_w+\mu^{-1}\K_{\bar
z}]=0.
\end{equation}
holds for all $\mu\in \C\setminus \{0\}$.  

If we assume the ASDYM connection $A$ is independent of $x_4$, and set $x=x_1$, $x_2=y$, and $x_3=t$,
then $A_w=\half(A_t-i\Phi)$ and $A_{\bar w}=\half(A_t+i\Phi)$,
where $\Phi=A_{x_4}$ is the Higgs field, $A=A_t\ dt + A_z\ dz +
A_{\bar z}\ d\bar z$ is a connection $1$-form on $\R^{2,1}$.  Then $(A,\Phi)$
satisfies the {\it space-time monopole equation\/}
$$
D_A\Phi= \ast F_A,
$$
where $\ast$ is the Hodge star operator with respect to the metric
$dx^2+dy^2- dt^2$. It has a  Lax pair induced from  \eqref{am}:
\begin{equation}\label{an}
\left[\half\K_t -\frac{i\phi}{2}
+\mu \K_z, \ \ \half\K_t +\frac{i\phi}{2} +\mu^{-1}\K_{\bar
z}\right]=0.
\end{equation}
Set
\begin{align*}
&D_1(\mu)= \half\K_t -\frac{i\phi}{2}
+\mu \K_z, \quad  D_2(\mu)= \half\K_t +\frac{i\phi}{2} +\mu^{-1}\K_{\bar z},\\
&\bca P_1(\mu) = D_1(\mu)-D_2(\mu)
= \frac{\mu-\mu^{-1}}{2}\K_x -i \frac{\mu+\mu^{-1}}{2}\K_y -i\phi,\\
P_2(\mu)= D_1(\mu)+ D_2(\mu)= \K_t +\mu \K_z + \mu^{-1}\K_{\bar z}.\eca
\end{align*} 
So \eqref{an} is equivalent to  
\begin{equation}\label{bl}
[\frac{\mu-\mu^{-1}}{2}\K_x - \frac{i(\mu+\mu^{-1})}{2}\K_y -i\phi,\   \K_t +\mu \K_z + \mu^{-1}\K_{\bar z}]=0.
\end{equation}
Note that the first operator is a linear operator in space variables. 
This is the Lax pair we use to construct monopoles with continuous scattering data. 

We need an equivalent form of the Lax pair to construct soliton monopoles.  First  we make a change of coordinates and spectral parameter: 
$$\xi=\frac{t+x}{2}, \quad \eta= \frac{t-x}{2}, \quad \mu=\frac{\tau-i}{\tau+i}.$$
  A direct computation shows that
\begin{align*}
L_1&=(\tau+i)D_1(\mu) + (\tau-i)D_2(\mu)=\tau\K_\xi -\K_y +\phi,\\
L_2&=\frac{1}{i}( (\tau+i)D_1(\mu) - (\tau-i)D_2(\mu))=\tau(\K_y+\phi) -\K_\eta.
\end{align*}
so $[\tau\K_\xi -\K_y +\phi, \ \tau(\K_y+\phi)-\K_\eta]=0$.  
Change spectral parameter again by $\l=\tau^{-1}$ to get
\begin{equation}\label{bd}
[\l(\K_y-\phi)-\K_\xi, \ \l\K_\eta-\K_y -\phi]=0.
\end{equation}
This is the Lax pair we use to construct B\"acklund transformations and  solitons for the monopole equation.  So we have

\bprop 
The following statements are equivalent for $(A,\phi)$:
\ben
\item $(A,\phi)$ is a solution of the space-time monopole equation,
\item \eqref{an} holds for all $\mu\in \C\setminus \{0\}$,
\item \eqref{bl} holds $\mu\in \C\setminus \{0\}$,
\item \eqref{bd} holds for all $\l\in \C$,
\item there exists $E_\mu(x,y,t)$ such that 
\begin{equation}\label{by}
\bca
(\frac{\mu-\mu^{-1}}{2}\,\p_x -\frac{i(\mu+\mu^{-1})}{2}\,\p_y)E_\mu=(\frac{\mu-\mu^{-1}}{2}\, A_x -\frac{i(\mu+\mu^{-1})}{2}\, A_y +i\phi)E_\mu,\\
(\p_t+\mu\p_z +\mu^{-1}\p_{\bar z})E_\mu= (A_t+ \mu A_z +\mu^{-1}A_{\bar z}) E_\mu,\\
E_{\bar\mu^{-1}}^* E_\mu=\I,
\eca
\end{equation}
\item there exists $\psi_\l(x,y,t)$ such that 
\begin{equation}\label{bz}
\bca
(\l\p_y -\p_\xi)\psi_\l = (\l(A_y+\phi) -A_\xi)\psi_\l,\\
(\l\p_\eta -\p_y)\psi_\l= (\l A_\eta - A_y+\phi)\psi_\l,\\
\psi_{\bar\l}^*\psi_\l=\I.
\eca
\end{equation}
\item
$E_\mu(x,y,t)$ is a solution of \eqref{by} if and only if 
$$\psi_\l(x,y,t)= E_{\frac{1-i\l}{1+i\l}}(x,y,t)$$ is a solution of \eqref{bz}.
\een
\eprop

We call solutions of \eqref{by} and \eqref{bz} {\it frames\/} of the monopole $(A,\phi)$.  
But frames are not unique.  In fact,  if $\psi_\l$ is a a solution of \eqref{bz} and $\phi_\l$ satisfies 
\begin{equation} \label{ca}
(\l\p_y - \p_\xi)\phi_\l =(\l\p_\eta-\p_y)\phi_\l=0, \quad \phi_{\bar \l}^*\phi_\l=\I,
\end{equation}
then $\psi_\l\phi_\l$ is also a solution of \eqref{bz}.  Moreover, given any meromorphic map $h:\C\to GL(n,\C)$ that satisfies $h(\bar\l)^*h(\l)=\I$, then $\phi_\l(x,y,t)= h(y+\l\xi +\l^{-1}\eta)$ is a solution of \eqref{ca}. However,  if $(A,\phi)$ is rapidly decaying in spatial variables, then we can choose normalizations (boundary conditions at infinity) so that there is a unique frame satisfying the normalization.     

\bs
\section{Direct scattering for flows in the $SU(n)$-hierarchy}

Let $V_n=\{ (\xi_{ij})\in su(n)\n \xi_{ii}=0 \ \forall \  1\leq i\leq n\}$. 
The phase space of evolution equations in the $SU(n)$-hierarchy is the set $\cs(\R, V_n)$ of all smooth $u:\R\to V_n$ that are rapidly decaying.  

Recall that $u$ is a solution of the $(b,j)$-flow in the $SU(n)$-hierarchy if and only if 
\begin{equation}\label{cc}
\bca 
\psi_\l^{-1}\p_x \psi_\l= a\l+ u,\\
\psi_\l^{-1}\p_t \psi_\l= b\l^j + Q_{b,1}(u)\l^{j-1} + \cdots + Q_{b,j}(u),\\
\psi_{\bar\l}^* \psi_\l=\I, 
\eca
\end{equation}
is solvable.  Since $u$ decays in $x$, it is natural to study solutions of the first linear operator in \eqref{cc} of the form $e^{a\l x} m(x,\l)$.   The direct scattering refers to the study of singularities of $m(x,\l)$ in spectral parameter $\l$. This was done by Beals and Coifman:

\bthm \label{ce}\cite{BeaCoi84} If $u\in \cs(\R, V_n)$, then there exist a bounded discrete subset $\D_u$ of $\C\setminus \R$ and a smooth map $m:\R\times \C\setminus (\R\cup \D_u)\to GL(n,\C)$ such that 
\ben
\item $\psi(x,\l)= e^{a\l x} m(x,\l)$ satisfies 
$d_x\psi= \psi(a\l+ u)$,
\item[(i)] $m(x,\bar\l)^*m(x,\l)=\I$ and $\lim_{x\to -\infty} m(x,\l)=\I$,
\item[(ii)]  $m(x,\l)$ is holomorphic for $\l\in \C\setminus (\R\cup \D_u)$, has poles at points in $\D_u$, and  $m_\pm(x, r)=\lim_{s\to 0^\pm} m(x,r+is)$ is smooth,
\item[(iii)] $m$ has an asymptotic expansion at $\l=\infty$:
$$m(x,\l) \sim \I + m_1(x)\l^{-1} + m_2(x)\l^{-2} + \cdots.$$
\een
Moreover,
\ben
\item there is an open dense subset $\cs_0(\R,V_n)$ of $\cs(\R,V_n)$ such that $\D_u$ is a finite set for  $u\in \cs_0(\R, V_n)$,  
\item if the $L^1$-norm of $u$ is less than $1$, then $m(x,\l)$ is holomorphic in $\l\in \C\setminus \R$, i.e., $\D_u$ is empty,
\item set $S(x,r)= m_+(x,r)m_-(x,r)^{-1}$, then $S^*=S$ and $S(x,r)-\I$ is rapidly decaying in $r$. 
\een
\ethm

The function $m$ in the above Theorem is called the {\it reduced wave function for the operator $d_x+ a\l + u$\/}, the poles and residues of $m$ are called the {\it discrete scattering data\/}, and the jump $S$ is called the {\it continuous scattering data\/} of $d_x+a\l+u$.

\bthm\cite{BeaCoi84, BeaCoi85}
Let $u$ be a solution of the $(b,j)$-flow \eqref{co} in the $SU(n)$-hierarchy such that $u(\cdot, t)\in \cs(\R, V_n)$, $m(\cdot, t, \cdot)$ and $S(\cdot,t,\cdot)$ the reduced wave function and  the continuous scattering data for $d_x+a\l + u(\cdot, t)$ respectively for each $t$. Set  
$\psi(x,t,\l)=e^{a\l x + b\l^j t} m(x,t,\l)$.  Then:
\ben
\item $\psi$ is a solution of \eqref{cc},
\item  $\bca \p_x S= [S, ar],\\ \p_t S= [S, br^j].\eca$ \hfil\break
In particular, 
$S(x,t,r)= e^{-(arx+br^j t)} s_0(r) e^{arx+br^j t}$
for some $s_0:\R\to GL(n,\C)$ such that $s_0^*=s_0$ and $s_0-\I$ is rapidly decaying.
\item If $u(\cdot, 0)$ has only continuous scattering data, then so is $u(\cdot, t)$. 
\item If the reduced wave function $m(\cdot, 0, \l)$ has a pole at $\l=\a$, then so is $m(\cdot, t, \l)$ for all $t$.
\item $u=[a, m_1]$, where $m_1$ is the coefficient of $\l^{-1}$ in 
the asymptotic expansion of $m(\cdot, \cdot,\l)$ at $\l=\infty$.  
\een
\ethm 

\bs
\section{Direct scattering for the space-time monopole equation}

The linear system associated to the Lax pair \eqref{bl} for the monopole equation is
\eqref{by}.  The first operator $P_1(\mu)$ is a linear operator in spatial variables only.  Given a rapidly decaying initial data $(A,\phi)$ on $\R^2$, the scattering data for the operator $P_1(\mu)$ is the singularity data of the solution $E_\mu$ for $P_1(\mu)E_\mu=0$ satisfying certain boundary condition.  
  
\bdefn A rapidly decaying spatial pair $(A,\phi):\R^2\to \oplus_{i=1}^4 su(n)$ is said to have only
{\it continuous scattering data\/} if there exists $E_\mu:\R^2\to GL(n,\C)$ defined on 
$\co^{\pm}_\e=\{\mu\in \C\n 1<|\mu|^{\pm 1}<1+\e\}$ for some $\e>0$ such that
\ben
\item
$
\bca P_1(\mu)E_\mu=\left(\frac{\mu-\mu^{-1}}{2}\K_x -
\frac{i(\mu+\mu^{-1})}{2}\K_y -i\phi\right)E_\mu =0,\\
E_\mu(\infty)=\I, \quad E_{{\bar u}^{-1}}=
(E_\mu^*)^{-1},\eca$
\item $\mu\mapsto E_\mu(x,t)$ are holomorphic in $\mu\in \co_\e^{\pm}$,
\item the limits $
\lim_{\mu\in \co^\pm, \mu\to e^{i\o}} E_\mu =
S_\o^\pm$ exist. 
\een
It follows from the reality condition that
$S_\o^-= {(S_\o^+)^*}^{-1}$. We call the non-negative Hermitian
matrix
$$S_\o=(S_\o^{-})^{-1}S_\o^+ = (S_\o^+)^\ast S_\o^+$$
the
{\it scattering matrix\/} or the {\it continuous scattering data\/}.
\edefn

Let $W^{2,1}$ denote the  space of maps $f$ whose partial derivatives up to second order are in $L^1$.

\bthm \label{bo} \cite{Vi90, FoIo01, DaTeUh06}
Assume that $(A,\phi)$ is a rapidly
decaying spatial data and $(A,\phi)$ is
small in $W^{2,1}$. Then the continuous scattering matrix $S_\o$ exists, $\I-S_\o$ decays for
each $\o$, and the scattering matrix $S_\o$ satisfies
\begin{itemize}
\item[(a)] $\I-S_\o$ is small in
$L^\infty$,
\item[(b)] $S^\ast_\o = S_\o \geq 0$,\item[(c)]
$ (-\sin\o \px +\cos\o
\py)S_\o=0$.
\end{itemize}
\ethm

\bthm \label{bn}\cite{Vi90, FoIo01, DaTeUh06}
If $(A,\phi)$ is a smooth solution of the space-time monopole equation in $\R^2\times (T_1,
T_2)$ and decays in spatial variables, and has a smooth continuous scattering
data. Then
$$
0= \left(\pt +\cos\o \px +\sin \o \py
\right)S_\o.
$$
Moreover, two gauge equivalent solutions give rise to the same scattering data.
\ethm

\bcor 
Let $(A,\phi)$ be as in Theorem \ref{bn}. Then there is a unique $s_0:\R\times S^1\to GL(n,\C)$ such that $s_0^*=s_0$, $s_0(r, e^{i\o})$ is rapidly decaying in $r\in \R$, and the continuous scattering data for $(A(\cdot, \cdot, t), \phi(\cdot, \cdot, t))$  is
$$S_\o(x, y, t)= s_0(x\cos\o + y\sin\o -t, e^{i\o}).$$
\ecor

\bs

\section{Inverse scattering for the $SU(n)$-hierarchy via loop group factorizations}

Given $u\in \cs(\R,V_n)$, the scattering data for the operator $L_u=d_x+a\l + u$ is the singularities of the reduced wave function, which contains two parts, the continuous (jumping line) and the discrete (poles) scattering data. The inverse scattering, which constructs $u$ from the scattering data of $L_u$, was done in \cite{ZakSha79, BeaCoi84}. 

By Theorem \ref{ce}, the scattering data only depends on $f(\l)=m(0,0,\l)$, where $m$ is the reduced wave function.   We identify the image of the scattering transform for those $u$'s with only continuous scattering data as a homogeneous space, and then use Pressley-Segal loop group factorization to construct the inverse scattering transform (\cite{TeUh98}).  

Let $\cd_-$ denote the group of smooth $f: \R\to GL(n,\C)$ such that
\ben
\item[(i)] $f$ is the boundary value of a holomorphic map in $\C_+=\{\l\in\C\n \Im(\l)>0\}$,
\item[(ii)] $f$ has the same asymptotic expansion at $r=\pm\infty$,
\item[(iii)] decompose $f(r)= p(r)v(r)$ with $p(r)$ upper triangular and $v(r)$ unitary, then $p-\I$ is rapidly decaying.
\een

Suppose $u$ is a solution of the $(b,j)$-flow \eqref{co} in the $SU(n)$-hierarchy with only continuous scattering data, and $m(x,t,\l)$ the reduced wave function for $u$.  Then $m(x,t,\cdot)\in \cd_-$.  
Set $f(\l)= m(0,0, \l)$.  Since $\psi(x,t,\l)= e^{a\l x+ b\l^j t} m(x,t,\l)$ satisfies \eqref{cc}, 
\begin{equation}\label{bg}
E(x,t,\l)= f(\l)^{-1} e^{(a\l x+ a\l^2 t)} m(x,t, \l),
\end{equation}
is a solution of \eqref{cc} with $E(0,0,\l)=\I$.    Because the right hand side of \eqref{bg} is holomorphic in $\C_+$,  $E(x,t,\l)$ is holomorphic in $\l\in \C_+$.  Proposition \ref{ci} implies that $E$ satisfies the $U(n)$-reality condition  $E(x,t,\bar\l)^*E(x,t,\l)=\I$. So by the reflection principal $E(x,t,\l)$ is holomorphic for all $\l\in \C$.

Set $e_{a,1}(x)(\l)= e^{a\l x}$, and $e_{b,j}(t)= e^{b\l^j t}$, $E(x,t)(\l)=E(x,t,\l)$, and $m(x,t)(\l)= m(x,t,\l)$.  Then
we can rewrite \eqref{bg} as
\begin{equation}\label{bh}
f^{-1}e_{a,1}(x) e_{b,j}(t)= E(x,t) m(x,t)^{-1}.
\end{equation}
Here $f, m(x,t)\in \cd_-$, and $e_{a,1}(x) e_{b,j}(t)$ and $E(x,t)$ holomorphic in $\C$ and satisfy the $U(n)$-reality condition. To construct the inverse scattering is to construct $m$ from $f$. In other words,   
given $f\in \cd_-$, we want to find a method to factor $f^{-1}e_{a,1}(x) e_{b,j}(t)$ as $E(x,t)m^{-1}(x,t)$ such that $E(x,t)$ satisfies the $U(n)$-reality condition and is holomorphic in $\C$ and $m(x,t)\in \cd_-$ for all $(x,t)$.  We need Pressley-Segal loop group factorization \cite{PrSe86} given below to do this factorization. 

Let $S^2\setminus S^1=\C\cup \{\infty\}= \W_+\cup \W_-$, where $\W_+=\{\mu\in \C\n |\mu|<1\}$ and $\W_-=\{\l\in S^2\n |\mu|>1\}$.     Let $\L(SL(n,\C))$ denote the group of smooth loops $g:S^1\to SL(n,\C)$, and $\L_+(SL(n,\C))$ the subgroup of $g\in \L(SL(n,\C))$ such that $g$ can be extended to a holomorphic map on $\W_+$ and $g(-1)$ is upper triangular with real diagonal entries.   Let $\L(SU(n))$ denote the loops in $SU(n)$.  The Pressely-Segal  factorization is the analogue of the  Iwasawa decomposition of $SL(n,\C)$ for loop groups:

\bthm[Pressely-Segal Factorization Theorem \cite{PrSe86}]\label{bq}
The multiplication map from $\L(SU(n))\times \L_+(SL(n,\C))$ to $\L(SL(n,\C))$ is a bijection. In particular, given $f\in \L(SL(n,\C))$, there exist unique $g\in \L(SU(n))$ and $h_+\in \L_+(SL(n,\C))$ such that $f=g h_+$.  
\ethm  

If we change the spectral parameter $\l$ by the linear fractional transformation $\mu= \frac{1+i\l}{1-i\l}$, then we can see that $\cd_-$ is isomorphic to a subgroup of $\L_+(SL(n,\C))$. In fact, we have

\bprop \cite{TeUh98}\label{br}  Given a map $g:S^1\to GL(n,\C)$, let  $\Phi(g):\R\to GL(n,\C)$ be the map defined by $\Phi(g)(r)=g(\frac{1+ir}{1-ir})$. Then: 
\ben
\item $g$ is smooth if and only if $\Phi(g)$ is smooth and has the same asymptotic expansion at $r=\pm\infty$.  
\item $j_\infty(g-\I)_{-1}=0$ (the infinite jet of $g-\I$ at $\mu=-1$) if and only if $\Phi(g)-\I$ is rapidly decaying.
\item Suppose $g$ extends  holomorphically to $|\mu|< 1$, and define $g$ on $|\mu|>1$ by $g(\mu)= (g(\bar\mu^{-1})^*)^{-1}$. Then  $f(\l)= g(\frac{1+i\l}{1-i\l})$ is holomorphic in $\l\in \C\setminus\R$ and satisfies the reality condition $f(\bar\l)^* f(\l)=\I$.  
\een
\eprop
 
 \bcor
$\cd_-$ is isomorphic to the subgroup of $g\in \L_+(SL(n,\C))$ such that $j_\infty(h-\I)_{-1}=0$ where $g= hv$ with $h$ upper triangular and $v$ unitary.  
\ecor

Now we go back to the problem of factorizing $f^{-1} e_{a,1}(x)$.  By Proposition \ref{br},   $\Phi^{-1}(f^{-1}e_{a,1}(x))$ does not belong to $\L(SL(n,\C))$.  So we can not use Theorem \ref{bq} to do the factorization directly.  However, if we write $f= pv$ with $p$ upper triangular and $v$ unitary, then by definition of $\cd_-$, $p-\I$ is rapidly decaying.  This implies that $\Phi^{-1}(e_{a,1}^{-1}(x)p^{-1}e_{a,1}(x))$ lies in $\L(SL(n,\C))$.  Apply the Pressley-Segal loop group factorization to get 
$$e_{a,1}^{-1}(x) p^{-1} e_{a,1}(x) = B(x) m(x)^{-1}$$ such that $\Phi^{-1}(B(x))\in \L(SU(n))$ and $\Phi^{-1}(m(x))\in \L_+(SL(n,\C))$.  
Since $p(\l)$, $e_{a,1}(x)(\l)$, and $m(x)(\l)$ are smooth for  $\l\in\R$ and can be extended holomorphically to $\l\in \C_+$, so is $B(x)(\l)$.  But $B(x)(\l)$ is unitary for $\l\in \R$ implies that $B(x)(\l)$ can be extended holomorphically across the real axis in the $\l$-plane by defining $B(x)(\l)= (B(x)(\bar\l)^*)^{-1}$.  Hence $\l\mapsto B(x)(\l)$ is holomorphic for all $\l\in \C$.  Therefore
\begin{align*}
f^{-1} e_{a,1}(x)&= v^{-1}p^{-1} e_{a,1}(x) = v^{-1} e_{a,1}(x) (e_{a,1}^{-1}(x) p^{-1} e_{a,1}(x))\\ &= v^{-1} e_{a,1}(x) B(x) m(x)^{-1} = E(x) m(x)^{-1}.
\end{align*}
But $E(x)(\l) = v^{-1}(\l) e^{a\l x} B(x)(\l)$ is holomorphic for $\l\in \C$.  

Since $E(x,\l)= f^{-1}(\l) e^{a\l x} m(x,\l)$, 
$$E^{-1}\p_x E= m^{-1} \p_x m + m^{-1}a\l m.$$
Use the asymptotic expansion at $\l=\infty$ to conclude that $E^{-1}\p_x E$ must be a degree one polynomial in $\l$.  So if $m_1(x)$ is the coefficient of $\l^{-1}$ in the asymptotic expansion of $m(x,\l)$ at $\l=\infty$, then 
$$E^{-1}\p_x E= a\l + u_f, \quad {\rm where\ } u_f=[a,m_1].$$
Note that the scattering data of $d_x+a\l +u_f$ is $e^{-a\l x} f_+f_-^{-1} e^{a\l x}$.  
However, the map $\cf(f)= u_f$ is not one to one. In fact, $u_{f_1}= u_{f_2}$ if and only if there is $h\in \cd_-$ such that $h(r)$ is diagonal for all $r\in \R$.  These give a rough idea of how the following results are obtained.

\bthm\label{bu}\cite{TeUh98}
Assume $a, b$ are diagonal matrices in $su(n)$, and $a$ has distinct eigenvalues.
If $f\in \cd_-$, then there exist $E(x,t,\l)$ and $m(x,t,\l)$ such that 
\ben 
\item $f^{-1}e_{a,1}(x)e_{b, j}(t) = E(x,t,\cdot) m(x,t, \cdot)^{-1}$,
\item $E$ is holomorphic for $\l\in \C$,  $E(x,t,\bar\l)^*E(x,t,\l)=\I$, and $m(x,t,\cdot)\in \cd_-$,
\item $u_f=[a, m_1]$ is a solution of the $(b,j)$-flow equation \eqref{co} in the $SU(n)$-hierarchy, and $E$ is the frame for the Lax pair associated to $u$ with initial condition $E(0,\l)=\I$, where $m_1(x,t)$ is the coefficient of $\l^{-1}$ in the asymptotic expansion of $m(x,t,\l)$ at $\l=\infty$,
\item $u_f(x,t)$ is defined for all $(x,t)\in \R^2$ and is rapidly decaying in $x$ for each $t$,
\item if $f$ also satisfies the $SU(n)/SO(n)$-reality condition, then $u_f$ is a solution of the $(b,j)$-flow in the $SU(n)/SO(n)$-hierarchy.
\een
\ethm

\bthm\cite{TeUh98}
Let $\cs^c(\R, V_n)$ denote the space of all $u\in \cs(\R, V_n)$ such that $L_u= d_x+ a\l + u$ has only continuous scattering data, and $\cd_-(A)$ denote the subgroup of $f\in\cd_-$ such that $f(r)$ is diagonal for all $r\in \R$, and $\cf:\cs^c(\R, V_n)\to \cd_-/\cd_-(A)$ defined by   $\cf(u)=[m(0, \cdot)]$, where $m(x,\l)$ is the reduced wave function of $L_u= d_x+ a\l + u$.  Then $\cf$ is a bijection, and $\cf^{-1}([f])=[a, m_1]$, where $m_1$ is the coefficient of $\l^{-1}$ in the asymptotic expansion of $m$ at $\l=\infty$.  \ethm

\bthm\cite{TeUh98}
Let $L^\tau_+(SL(n,\C))$ denote the group of holomorphic maps $f:\C\to GL(n,\C)$ that satisfy the $SU(n)$-reality condition, and  $\cd_+(A)$  the subgroup of $L^\tau_+(SL(n,\C))$ generated by 
$\{e_{b,j}(t)\n b\in su(n) \ {\rm diagonal\/},\  j\geq 1 \ {\rm integer}\}$.
 Then $\cd_+(A)$ acts on $\cd_-/\cd_-(A)$ by $e_{a,j}(t)\ast [f]= [m(t)]$, where $m(t)$ is obtained by factoring $f^{-1}e_{b,j}(t)= E(t)m(t)^{-1}$ such that $E(t)\in L^\tau_+(SL(n,\C))$ and $m(t)\in \cd_-$.  Moreover, 
 the $(b, j)$-flow in the $SU(n)$-hierarchy corresponds to the action of $e_{b,j}(t)$ on $\cd_-/\cd_-(A)$ under the isomorphism $\cf$.
 \ethm
 
 \bthm \cite{TeUh98}
  Let $a_1, a_2, \ldots, a_{n}$ be linearly independent diagonal matrices in $u(n)$, and $f\in \cd_-$. Then we can factor 
  $$f^{-1}e_{a_1,1}(x_1) \cdots e_{a_{n},1}(x_{n}) = E(x) m(x)^{-1}$$ such that $E(x)\in L_+^\tau(GL(n,\C))$ and $m(x)\in \cd_-$. Moreover, 
  \ben
  \item $v=m_1^\perp$ is a solution of the $U(n)$-system, where $m_1(x)$ is the coefficient of $\l^{-1}$ in the asymptotic expansion of $m(x)(\l)$ at $\l=\infty$ and $\xi^\perp=\xi-\sum_{i=1}^n \xi_{ii} e_{ii}$,
  \item if $f\in \cd_-$ satisfies the $\frac{U(n)}{O(n)}$-reality condition, then $v=(m_1)^\perp$ is a solution of the $\frac{U(n)}{O(n)}$-system.  
  \een
\ethm

In other words, the $U(n)$-system is the system obtained by putting the $(a_1, 1)$-, $\ldots$, $(a_n,1)$-flow in the $U(n)$-hierarchy together.
  
  \bs
  \section{The inverse scattering for  monopole equations}

The scattering data of the linear operator 
$$P_1(\mu)= \mu\K_z -\mu^{-1}\K_{\bar z} -i\phi$$ on the $z=x+iy$ plane for the rapidly decaying spatial pair $(A,\phi)$ is a smooth map $s_0:\R^1\times S^1\to GL(n,\C)$.  The inverse scattering for $P_1(\mu)$, which constructs $(A, \phi):\R^2\to \oplus^4 su(n)$ from $s_0$,  was done in \cite{Vi90, FoIo01}.  In this section, we give a brief review of the construction of the inverse scattering transform for $P_1(\mu)$ via Pressley-Segal  loop group factorization  given in \cite{DaTeUh06}.

\bthm\label{cg}\cite{DaTeUh06}
Suppose $s:\R\times S^1\to GL(n,\C)$ is smooth such that $s^*=s\geq 0$ and $s(r, e^{i\o})-\I$ is rapidly decaying for $r\in \R$.  Define 
$$S(x,y,t,e^{i\o})= s(x\cos\o +y\sin\o -t, e^{i\o}).$$
Then there exists a smooth $E:\R^{2,1}\times (\C\setminus S^1)\to GL(n,\C)$ such that 
\ben
\item $E_{\bar\mu^{-1}}^*E_\mu=\I$, where $E_\mu= E(\cdots, \mu)$,
\item $((\p_t+\mu\p_z)E_\mu)E_\mu^{-1} = B_0+\mu B_1$ and $((\p_t +\mu^{-1}\p_{\bar z}) E_\mu)E_\mu^{-1}= -(B_0^*+\mu^{-1} B_1^*)$ for some $B_0, B_1:\R^{2,1}\to sl(n,\C)$, 
\item Set $A_z= \half B_1, \quad A_t =\half (B_0-B_0^*), \quad \phi=\frac{i}{2}(B_0+B_0^*)$, then $(A,\phi)$ is a solution of the space-time monopole equation decaying rapidly in the spatial variables,
\item the scattering data for $(A(\cdot, \cdot, t), \phi(\cdot, \cdot, t))$ is $S(\cdot,\cdot, t,\cdot)$.  
\een
\ethm

Here is a sketch of the proof: Set $S_\mu= S(\cdots, \mu)$ for $|\mu|=1$.  Write $S_\mu=P_\mu^2$ with $P_\mu^*=P_\mu$.  By Pressley-Segal factorization  Theorem \ref{bq} we can factor $P_\mu= U_\mu E_\mu^+$ with $U_\cdot$ a loop in $SU(n)$  and $E^+_\mu$ extends holomorphically to $|\mu|<1$.  Define $E^-_\mu = ((E_{\bar\mu^{-1}}^+ )^*)^{-1}$ for $|\mu|>1$.  Then $S_\mu= (E_\mu^-)^{-1} E_\mu^+$.  The rest of the Theorem can be proved using the fact that $(-\sin\o\ \p_x +\cos\o\ \p_y)S=0$ and $(\cos\o\ \p_x + \sin\o\ \p_y-\p_t)S=0$.  

\bcor \cite{DaTeUh06} Suppose $(A_0,\phi_0)$ is a rapidly decaying spatial pair with only continuous scattering data.  Then there is a global solution $(A, \phi)$ of the space-time monopole equation decaying rapidly in the spatial variables such that the scattering data of $(A(\cdot, \cdot, 0), \phi(\cdot, \cdot,0))$ and $(A_0, \phi_0)$ are the same.  Moreover, any two such solutions are gauge equivalent.
\ecor

\bcor  There is  a bijective correspondence between the space of solutions of the space-time monopole equation with only continuous scattering data modulo the gauge group, and the group of maps $f:\R\to \L^\tau_+(SL(n,\C))$ such that $(f^*f)(r)(e^{i\o})-\I$ is rapidly decaying in $r\in \R$. 
\ecor 

\bs 
\section{Birkhoff factorization and local solutions}

The factorization \eqref{bh} 
$$f^{-1}(\l) e^{a\l x + b\l^j t}= E(x,t,\l) m(x,t,\l)^{-1}$$
is for $f, m(x,t,\cdot)$ in $\cd_-$ and $E$ holomorphic in $\l\in \C$ and satisfies the $U(n)$-reality condition.  However if $f$ is holomorphic at $\l=\infty$, then we can use the Birkhoff factorization to get $E$ and $m$ such that $E$ is holomorphic in $\C$ and $m$ is holomorphic at $\l=\infty$.  Moreover, it can be shown easily that $E$ is the frame of a solution of the $(b,j)$-flow with $E(0,0,\l)=\I$.  Since the Birkhoff factorization only works on an open dense subset of loops, solutions constructed this way are local solutions defined in a neighborhood of $(0,0)$. 

   Let $\e>0$, and $\co_\e=\{\l\n |\l | > \frac{1}{ \e}\}$  an
open neighborhood of $\infty$ in $S^2=\C\cup\{\infty\}$. Then $S^2=C\cup
\co_{\e}$.  Let  
$L^\tau(SL(n,\C))$ denote the group of holomorphic maps $f$ from  $C\cap \co_{\infty}$
to $SL(n,\C)$ satisfying the $SU(n)$-reality condition $f(\bar\l)^*f(\l)= \I$, $L^\tau_+(SL(n,\C))$ the subgroup of $f\in L^\tau(SL(n,\C))$ that extend
holomorphically to $\C$, and $L^\tau_-(SL(n,\C))$ the subgroup of $f\in L^\tau(SL(n,\C))$ that
extend holomorphically to $\co_\e$ and satisfying $f(\infty)=I$.   

\bthm[\bf Birkhoff Factorization Theorem]  \label{bs} (cf. \cite{PrSe86})
The multiplication map $L^\tau_+(SL(n,\C))\times L^\tau_-(SL(n,\C))\to L^\tau(SL(n,\C))$ is injective and the image is an open dense subset of $L^\tau(SL(n,\C))$.
\ethm

Let $a, b$ be diagonal matrices in $su(n)$ such that $a$ is regular.  Then $e_{a,1}(x) e_{b,j}(t)\in L^\tau_+(SL(n,\C))$.   Given  $f\in L_-^\tau(SL(n,\C))$,  by the Birkhoff factorization there exists $\d>0$ such that
$$f^{-1}e_{a,1}(x) e_{b,j}(t) = E(x,t) m(x,t)^{-1}$$
with $E(x,t)\in L_+^\tau(SL(n,\C))$ and $m(x,t)\in L_-^\tau(SL(n,\C))$ for all $(x,t)\in B_\d(0)$. Here $B_\d(0)$ is the ball of radius $\d$ centered at $(0,0)$. Then 
$$\bca
E^{-1} \p_x E= m^{-1}m_x + m^{-1} a\l m, \\
E^{-1}\p_t E= m^{-1}m_x + m^{-1} b\l^j m.\eca$$
Since $m$ is holomorphic at $\l=\infty$ and $m(x,t)(\infty)=\I$, $E^{-1}\p_x E$ and $E^{-1}\p_t E$ must be a polynomial of degree $1$ and $j$ in $\l$ respectively.  Hence $E$ must be a frame of a solution of the $(b,j)$-flow in the $SU(n)$-hierarchy.  So we have 

\bthm \label{bp}\cite{TeUh98}
If $f\in L_-^\tau(SL(n,\C))$, then there exist an open neighborhood $\co$ of  $(0,0)$, $E(x,t)\in L^\tau_+(SL(n,\C))$, and $m(x,t)\in L^\tau_-(SL(n, \C))$ such that
$f^{-1}e_{a,1}(x) e_{b,j}(t)= E(x,t) m(x,t)^{-1}$ for all $(x,t)\in \co$.  Moreover,
\ben 
\item $u=[a, m_1]$ is a solution of the $(b,j)$-flow in the $SU(n)$-hierarchy, where $m_1(x,t)$ is the coefficient of $\l^{-1}$ the expansion of $m(x,t)(\l)$ at $\l=\infty$,  (we will use $f\ast 0$ to denote $u$),
\item $E$ is the frame of the Lax pair of $u$ such that $E(0)(\l)= \I$,
\item if $f$ satisfies the $SU(n)/SO(n)$-reality condition, then $u=[a,m_1]$ is a solution of the $(b,j)$-flow in the $SU(n)/SO(n)$-hierarchy.
\een
\ethm

Note that for $f\in \cd_-$, Theorem \ref{bu} gives a global solution $u_f$ of the $(b,j)$-flow on $\cs(\R, V_n)$ in the $SU(n)$-hierarchy.  But for $f\in L_-^\tau(SL(n,\C))$, the above theorem only gives a local solution of the $(b,j)$-flow in general. 

\bthm\cite{TeUh98}\label{bw}
If $f\in L_-^\tau(GL(n,\C))$, then there exist an open neighborhood $\co$ of $0$ in $\R^n$, $E(x)\in L^\tau_+(GL(n,\C))$ and $m(x)\in L^\tau_-(GL(n, \C))$ for $x\in \co$ such that
$f^{-1}e_{a_1,1}(x_1)\cdots e_{a_n, 1}(x_n)= E(x) m(x)^{-1}.$  Moreover,
\ben 
\item $v=m_1^\perp$ is a solution of the $U(n)$-system, where $m_1(x)$ is the coefficient of $\l^{-1}$ in the expansion of $m(x)(\l)$ at $\l=\infty$ and $m_1^\perp=m_1-\sum_i (m_1)_{ii} e_{ii}$,  (we will use $f\ast 0$ to denote $v$),
\item $E$ is the frame of the Lax pair \eqref{cf} of $v$ such that $E(0)(\l)=\I$,
\item If $f$ also satisfies the $U(n)/O(n)$-reality condition, then $v=f\ast 0$ is a solution of the $U(n)/O(n)$-system. 
\een
\ethm

\bs
\section{B\"acklund transformations for the $U(n)$-hierarchy}

In general, the solution $f\ast 0$ constructed in Theorem \ref{bp} has singularities.  But if $f$ is rational, then $f\ast 0$ is a global solution of the $(b,j)$-flow in the $SU(n)$-hierarchy, and can be computed explicitly.  These are the soliton solutions.  Moreover, if $f$ is rational with only one simple pole and $E$ is a frame of a solution $u$, then the Birkhoff factorization $fE = \ti E \ti f$ can be carried out by an explicit algebraic algorithm so that $\ti E$ is a frame of  the new solution.  This give B\"acklund transformations for the $(b,j)$-flow.  

If $\a\in \C\setminus \R$ and $\pi$ is a Hermitian projection  of $\C^n$, then the map
$$g_{\a, \pi}(\l)= \I + \frac{\a-\bar \a}{\l-\a}\ \pi$$
satisfies the $U(n)$-reality condition. So  $g_{\a, \pi}\in L_-^\tau(GL(n,\C))$. 
 
 The following  Theorem is a key ingredient for constructing B\"acklund transformations for the $(b,j)$-flow in the $SU(n)$-hierarchy.

\bthm\label{aw}\cite{TeUh00}
Given $f\in L_+^\tau(SL(n,\C))$ and $g_{\a,\pi}$, let $\ti\pi$ be the Hermitian projection of $\C^n$  onto $f(\a)^{-1}(\Im\pi)$.  Then $g_{\a, \pi} f= \ti f g_{\a,\ti\pi}$ and $\ti f\in L^\tau_+(SL(n,\C)$.
\ethm

\begin{proof} Set 
$\ti f(\l)= g_{\a,\pi}(\l) f(\l) g_{\a,\ti\pi}(\l)^{-1} = (\I+\frac{\a-\bar \a}{\l-\a}\pi^\perp)f(\l) (\I + \frac{\bar\a-\a}{\l-\bar \a}\ti\pi^\perp)$.
Note that $\ti f$ is holomorphic for $\l\in \C\setminus\{\a, \bar \a\}$. But
$${\rm Res}(\ti f, \a)=(\a-\bar\a)  \pi^\perp f(\a) \ti \pi, \quad {\rm Res}(\ti f, \bar\a)=(\bar\a-\a)\pi f(\bar\a)\ti\pi^\perp.$$
 By definition $f(\a)(\Im\ti\pi)=\Im \pi$, so $\ti f$ is holomorphic at $\l=\a$.  Set $V=\Im\pi$ and $\ti V=\Im\ti\pi$.   Since $f$ satisfies the reality condition,  we have
 $$(f(\bar\a)(\ti V^\perp), V)= (\ti V^\perp, f(\bar \a)^*(V))= (\ti V^\perp, f(\a)^{-1}(V)) =(\ti V^\perp, \ti V)=0.$$
 This implies that ${\rm Res\/}(\ti f,\bar \a)=0$, hence $\ti f$ is holomorphic in $\C$.  
\end{proof}

The proof of the above theorem in fact gives the following more general result:

\bthm\label{as} \cite{DaTe04}
Let $\co$ be an open subset of $\C$ that is invariant under complex conjugation,
 and $f:\co\to SL(n,\C)$ a meromorphic map satisfying the $U(n)$-reality condition.  Let $\a\in \C\setminus \R$, and $\pi$ a Hermitian projection of $\C^n$.  Suppose $f$ is holomorphic and non-singular at $\l=\a$. Let $\ti \pi$ denote the Hermitian projection of $\C^n$ onto $f(\a)^{-1}(\Im \pi)$.  Then
 $g_{\a, \pi} f= \ti f g_{\a, \ti \pi}$, and  $\ti f$ is holomorphic and non-degenerate at $\l=\a$ and satisfies the $U(n)$-reality condition.
\ethm

\bthm[\bf B\"acklund transformation for the $(b,j)$-flow]\cite{TeUh00} \label{bx}
Suppose $u$ is a solution of the $(b,j)$-flow \eqref{co} in the $SU(n)$-hierarchy, and $E(x,t,\l)$ is the frame for the Lax pair of $u$ such that $E$ is holomorphic for $\l\in \C$ and $E(0,0, \l)=\I$.  Given $\a\in \C\setminus \R$ and a Hermitian projection $\pi$ of $\C^n$, set $\ti \pi(x,t)$ to be the Hermitian projection of $E(x,t,\a)^{-1}(\Im\pi)$, $\ti E= g_{\a, \pi} E g_{\a, \ti\pi}^{-1}$, and $\ti u= u+ (\a-\bar\a)[a, \ti\pi]$. Then 
\ben
\item $\ti u$ is again a solution of \eqref{co}, and $\ti E$ is the frame of $\ti u$,
\item if $u$ is smooth for all $(x,t)\in \R^2$, then so is $\ti u$,
\item if $u(x,t)$ is rapidly decaying in $x$ for all $t$, then so is  $g_{\a,\pi}\ast u$.
\een
\ethm

Let $\cd_-^r$ denote the group of rational maps $g:S^2\to GL(n,\C)$ that satisfy the $U(n)$-reality condition and $g(\infty)=\I$.
Uhlenbeck proved in \cite{Uh89} that $\cd_-^r$ is generated by the set 
$\{g_{\a,\pi}\n \a\in \C\setminus \R, \ \pi \ {\rm is \ a \ Hermitian\ projection\  of \ } \C^n\}$.

A {\it pure soliton\/} for the $(b,j)$-flow is a solution that is rapidly decaying in the spatial variable, has no continuous scattering data and has finitely many discrete scattering data, so its reduced wave function $m(x,\l)$ is  rational in $\l$. Or equivalently, its reduced wave function $m$ lies in the group $\cd_-^r$.  

\bcor \cite{TeUh00}
The group $\cd_-^r$ acts on the space of solutions of the $(b,j)$-flow in the $SU(n)$-hierarchy.  In fact, if $g=g_{\a_1, \pi_1}\cdots g_{\a_k, \pi_k}$, then $g\ast u= g_{\a_1, \pi_1}\ast(\cdots \ast (g_{\a_k.\pi_k}\ast u)\cdots)$. 
\ecor

\bcor\cite{TeUh98}
Let $E_0(x,t,\l)= e^{a\l x + b\l^j t}$
 ($E_0$  is the frame for the vacuum solution $u=0$ of the $(b,j)$-flow).  If we apply BT (Theorem \ref{bx}) to $E_0$ repeatedly, then we obtain all pure soliton solutions of the $(b,j)$-flow, i.e., solutions with continuous scattering data $S=\I$ and finitely many discrete scattering data.
\ecor
  
  \bthm\cite{TeUh98}
Let $u$ be the global solution of the $(b,j)$-flow \eqref{co} in the $SU(n)$-hierarchy constructed in Theorem \ref{bu} with only continuous scattering data, and $E$ its frame.  If we apply BT (Theorem \ref{bx}) repeatedly to $E$, then we obtain solutions of the $(b,j)$-flow  that have both continuous and finite discrete scattering data.  Conversely, any solution $u$ of the $(b,j)$-flow in the $SU(n)$-hierarchy that has continuous scattering data and finite discrete scattering data can be constructed this way. 
\ethm

\bthm[\bf B\"acklund transformation for the $U(n)$-system]\cite{TeUh00}
Suppose $v$ is a solution of the $U(n)$-system, and $E(x,\l)$ is the frame for the Lax pair of $v$ such that $E$ is holomorphic for $\l\in \C$ and $E(0,\l)=\I$.  Given $\a\in \C\setminus \R$ and a Hermitian projection $\pi$ of $\C^n$, set $\ti \pi(x)$ to be the Hermitian projection of $E(x,\a)^{-1}(\Im\pi)$, $\ti E= g_{\a, \pi} E g_{\a, \ti\pi}^{-1}$, and $\ti v= v+ (\a-\bar\a)\ti\pi^\perp$, where $\xi^\perp=\xi-\sum_{i=1}^n \xi_{ii} e_{ii}$. Then 
\ben
\item $\ti v$ is again a solution of the $U(n)$-system, and $\ti E$ is a frame of the Lax pair of $\ti v$,
\item if $v$ is smooth for all $x\in \R^n$, then so is $\ti v$.
\een
\ethm

The group $\cd_-^r$ acts on the space of solutions of the $U(n)$-system such that $g_{\a, \pi}\ast v=\ti v$, where $\ti v$ is given in the above theorem.  

It is easy to see that if $u$ is a solution of the $\unon$-system, $s\in \R$, and $\bar\pi=\pi$, then $g_{is,\pi}$ satisfies the $U(n)/O(n)$-reality condition and $g_{is,\pi}\ast u$ is also a solution of the $\unon$-system.  In general,

\bcor \cite{TeUh00}
 If $g\in \cd_-^r$ satisfies the $\unon$-reality condition and $v$ is a solution of the $\unon$-system, then $g\ast v$ is again a solution of the $\unon$-system. 
\ecor

\bs
\section{B\"acklund transformations for the space-time monopole equation}

We use Lax pair \eqref{bd} to construct soliton solutions and B\"acklund transformations for the monopole equation.  Since the spectral parameter $\l$ in \eqref{bd} is related to the spectral parameter  $\mu$ in \eqref{bl} by $\mu=\frac{1-i\l}{1+i\l}$, the continuous scattering data for \eqref{bd} is the jump  across the real axis, and the discrete scattering data is given by the poles in $\C\setminus \R$ and their residues.   

\bdefn  A monopole $(A,\phi)$ rapidly decaying in the spatial variable is a {\it $k$-soliton\/} if there is a gauge equivalent monopole with a frame (solution of \eqref{bz})  $\psi_\l$ that is rational in $\l$ with $k$ poles, $\psi_{\infty}=\I$, and  $\lim_{||(x,y)||\to \infty} \psi_\l(x,y,t,\l)= h(\l)$ is independent of $t$.  
\edefn

We identify the set of all rank $k$ Hermitian projections of $\C^n$ as the complex Grassmannian $\Gr(k, \C^n)$ by $\pi\mapsto \Im\pi$.  

\bthm\cite{Wa88}  Let $\a\in \C\setminus \R$,  $\pi_0:S^2\to\Gr(k, \C^n)$  a holomorphic map, $\xi= \half(t+x)$, $\eta= \half (t-x)$, $\pi(x,y,t)= \pi_0(y+\a \xi + \a^{-1} \eta)$, and
$$g_{\a, \pi_0}(x,y,t)= \I + \frac{\a-\bar \a}{\l-\a}\ \pi(x,y,t).$$
Then $g_{\a, \pi}$ is a $1$-soliton monopole frame.  Moreover, all  $1$-soliton frames are of this form up to gauge equivalence. 
\ethm

The following theorem is a consequence of Theorem \ref{as} and the fact that  $\l\p_y-\p_\xi$ and $\l\p_\eta-\p_y$ are derivations. 

\bthm [{\bf BT for Monopoles}]\label{bc}\cite{DaTe04}
Suppose $\a\in \C\setminus \R$ is a constant, and $\psi$ is a frame of the monopole solution $(A,\phi)$ (i.e.,  solution of \eqref{bz}), and $\psi(x,y,t,\tau)$ is holomorphic and non-degenerate at $\tau= \a$.  Let $g_{\a, \pi}$ be a $1$-soliton monopole frame, $\ti\pi(x,y,t)$ the Hermitian projection of $\C^n$ onto 
$\psi(x,y,t,\a)(\Im\pi(x,y,t))$, and $\ti\psi= g_{\a, \ti\pi} \psi g_{\a, \pi}^{-1}$.
Then 
\ben
\item $\ti\psi$ is holomorphic and non-degenerate at $\tau=\a$,
\item $\psi_1= g_{\a, \ti\pi} \psi= \ti\psi g_{\a,\pi}$ is a frame for \eqref{bz} with $\ti A, \ti \phi$ given by 
$$\bca
\ti A_\eta= A_\eta,\\
\ti A_\xi = (1-\frac{\bar a}{\a}) (\p_\xi\ti \pi) h + h^{-1}A_\xi h,\\
\ti A_y +\ti\phi= A_y+\phi,\\
\ti A_y -\ti \phi= (1-\frac{\bar\a}{\a})(\p_y\ti\pi) h + h^{-1}(A_y-\phi) h,
\eca
$$
where $h=\ti\pi +\frac{\a}{\bar\a} \ \ti\pi^\perp$, 
\item $(\ti A, \ti \phi)$ is a solution of the space-time monopole equation.
\een
\ethm

If we apply Theorem \ref{bc} to a $1$-soliton $k$-times, then we get a $(k+1)$-soliton whose frame has $(k+1)$ distinct simple poles.  Moreover, we have

\bcor
Suppose  $(A,\phi)$ is a solution of the space-time monopole equation with only continuous scattering data and $E$ is its frame constructed in Theorem \ref{cg}. If we apply Theorem \ref{bc} to $E$ repeatedly, then we obtain a monopole whose frame has both continuous  scattering data and finitely many distinct simple poles. 
\ecor

Note that a BT for flows in the $SU(n)$-hierarchy  adds to a given solution, a soliton with scattering pole at $\a$, regardless of whether the given solution already has a scattering pole at $\a$ or not.  But this is not the case for the monopole equation, so BTs produce soliton monopole frames with distinct simple poles only.   Ward and his group (\cite{Wa95, Io96, An97, IoZa98a}) take limits of soliton monopole frames with $2$ and $3$ distinct poles to construct $2$- and $3$- soliton monopoles with a double  and a triple pole at $i$ that are time dependent. 
 Dai and Terng used BT (Theorem \ref{bc}) and a systematic limiting method  to construct rational monopole frames with arbitrary poles and multiplicities:

\bthm\cite{DaTe04}  Given $\a_i\in \C\setminus \R$ and positive integers $n_i$ for $1\leq i\leq k$, there are soliton monopole frames that are rational and have poles at $\a_1, \ldots, \a_k$ with multiplicities $n_1, \ldots, n_k$.
\ethm

Below is a more general B\"acklund transformation that adds a multiplicity $k$ pole at $\l=\a$ to a given monopole frame.
    
  \bthm \cite{DaTe04}  Suppose $\psi$ is a monopole frame that is holomorphic and non-degenerate at $\l=\a$ and $\phi$ is a soliton (rational) monopole frame with a single pole at $\l=\a$ with multiplicity $k$.  Then there exist unique $\ti\phi, \ti \psi$ such that $\psi_1= \ti\phi \psi = \ti \psi \phi$ is a monopole frame, $\ti \phi$ is rational with a single pole at $\l=\a$ with multiplicity $k$, and $\ti\psi$ is holomorphic and non-degenerate at $\l=\a$.    We use $\phi \ast \psi$ to denote $\psi_1$.  
   \ethm
   
 \bthm\cite{DaTeUh06} If $\psi$ is a monopole frame with both continuous scattering data and finitely many poles, then there exist unique monopole frames $\psi_c$ and $\phi$ such that $\psi_c$ has only continuous scattering data, $\phi$ has only discrete scattering data, and $\psi= \phi\ast \psi_c$.  
 \ethm

\bs


\begin{thebibliography}{99}

\bibitem{AKNS74}
Ablowitz, M.J., Kaup, D.J., Newell, A.C. and Segur, H., \emph{{T}he inverse
scattering transform - Fourier analysis for nonlinear
problems}, Stud. Appl. Math. \textbf{53} (1974), 249--315

\bibitem{AblCla91}
Ablowitz, M.J.,Clarkson, P.A.,\emph{{S}olitons,
non-linear evolution equations and inverse scattering},  Cambridge Univ.
Press (1991)

\bibitem{An97}
Anand, C. K., \emph{{W}ard's solitons},  Geom. Topol., {\textbf 1}
(1997), 9--20.

\bibitem{BeaCoi84}
Beals, R., Coifman, R.R.,\emph{{S}cattering and inverse scattering for
first order systems}, Commun. Pure Appl. Math. \textbf{37} (1984), 39--90

\bibitem{BeaCoi85} 
Beals, R., Coifman, R.R., \emph{{I}nverse scattering and evolution
equations}, Commun. Pure Appl. Math., \textbf{38} (1985), 29-42

\bibitem{BeaCoi85b} 
Beals, R., Coifman, R.R., \emph{{M}ultidimensional inverse scattering and nonlinear partial differential equations}, Proc. Symp. Pure Math., \textbf{43} (1985), 45-70

\bibitem{BeaCoi89}
Beals, R., Coifman, R.R., \emph{{L}inear spectral problems, non-linear
equations and the $\bar\partial$-method}, Inverse Problems, \textbf{5} (1989), 87-130

\bibitem{BelZak78}
Belavin, A. A., Zakharov, V.E., \emph{{Y}ang-Mills equations as an inverse scattering problem}, Phys. Lett.. \textbf{73B} (1978), 53-57

\bibitem{BrDuPaTe02}
Br\"uck, M., Du, X., Park, J., and Terng, C.L., \emph{{S}ubmanifold geometry of real 
Grassmannian systems}, The Memoirs, vol 155, AMS,  \textbf{735} (2002), 1--95

\bibitem{Bur00}
Burstall, F., \emph{{I}sothermic surfaces: conformal geometry, Clifford
algebras and integrable systems}, to appear in ``Integrable systems, Geometry, and
Topology'',  AMS-International Press, math-dg/0003096

\bibitem{CieGolSym95}
Cie\'sli\'nski, J., Goldstein, P.,  and Sym, A., \emph{{I}sothermic surfaces in
  {${E}^3$} as soliton surfaces}, Phys.\ Lett.\ A \textbf{205} (1995), 37--43
  
\bibitem{DaTe04}
Dai, B. and Terng, C. L., \emph{{B}\"acklund transformation, Ward
solitons, and unitons}, to appear in J. Differential Geometry,  arXiv:math.DG/0405363.

\bibitem{DaTeUh06}
Dai, B., Terng, C. L., and Uhlenbeck, K., \emph{{O}n the space-time monopole equations}, to appear in J. Differential Geometry Survey, arXiv:math.DG/0602607

\bibitem{FadTak87}
Faddeev, L.D., Takhtajan, L.A., \emph{{H}amiltonian methods in the theory of solitons}, (1987),
 Springer-Verlag

\bibitem{FerPed96b}
Ferus, D., Pedit, F., \emph{{I}sometric immersions of space forms and
soliton theory}, Math. Ann., \textbf{305} (1996), 329--342

\bibitem{FoIo01}
Fokas, A. S. and Ioannidou, T. A., \emph{{T}he inverse spectral
theory for the Ward equation and for the $2+1$ chiral model},
Comm. Appl. Anal., {\textbf 5} (2001), no. 2, 235--246.

\bibitem{GGKM67}
Gardner, C.S., Greene, J.M., Kruskal, M.D., Miura, R.M., \emph{{M}ethod for
solving the Korteweg-de Vries equation}, Physics Rev. 
Lett. \textbf{19} (1967), 1095-1097

\bibitem{HarSaiShn84b} 
Harnad, J., Saint-Aubin, Y., Shnider, S., \emph{{T}he Soliton Correlation Matrix and  the Reduction Problem for Integrable Systems}, Commun. Math. Phys. \textbf{93} (1984), 33-56

\bibitem{Io96}
Ioannidou, T., \emph{{S}oliton solutions and nontrivial scattering
in an integrable chiral model in $(2+1)$ dimensions}, J. Math.
Phys., {\textbf 37} (1996), 3422--3441.

\bibitem{IoZa98a}
Ioannidou, T. and Zakrzewski, W., \emph{{S}olutions of the
modified chiral model in $(2+1)$ dimensions}, J. Math. Phys.,
{\textbf 39} (1998) no.5, 2693--2701.

\bibitem{MaZa81}
Manakov, S. V. and Zakharov, V. E., \emph{{T}hree-dimensional
model of relativistic-invariant theory, integrable by the inverse
scattering transform}, Lett. Math. Phys., {\textbf 5} (1981),
247--253.

\bibitem{PrSe86}
Pressley, A. and Segal, G., \emph{{L}oop groups}, Oxford
University Press, 1986.

\bibitem{Sat84}
Sattinger, D.H., \emph{{H}amiltonian hierarchies on semi-simple Lie
algebras}, Stud. Appl. Math., \textbf{72} (1984), 65--86

\bibitem{Ten98} 
Tenenblat, K., \emph{{T}ransformations of manifolds and
applications to differential equations}, Pitman Monographs and Surveys
in Pure and Applied Mathematics, \textbf{93} (1998), Longman, Harlow


\bibitem{Ter97}
Terng, C.L., \emph{{S}oliton equations and differential
geometry}, J. Differential Geometry, \textbf{45} (1997), 407--445

\bibitem{Ter03}
Terng, C.L., \emph{Geometries and symmetries of soliton equations and integrable elliptic systems}, to appear in 
Surveys on Geometry and Integrable Systems,  Advanced Studies in Pure Mathematics, Mathematical Society
of Japan, math.DG/0212372

\bibitem{TeUh98}
Terng, C. L. and Uhlenbeck, K.,  \emph{{P}oisson actions and
scattering theory for integrable systems}. Surveys in differential
geometry: integrable systems, 315--402, Surv. Diff. Geom., IV,
International Press, Boston, MA, 1998.

\bibitem{TeUh99}
Terng, C.L., Uhlenbeck, K., \emph{{S}chr\"odinger flows on Grassmannians}, Integrable systems, Geometry, and Topology, AMS/IP Stud. Adv. Math. 36 (2006), 235-256, math.DG/9901086

\bibitem{TeUh00}
Terng, C. L. and Uhlenbeck, K.,  \emph{{B}\"{a}cklund
transformations and loop group actions}, Comm. Pure Appl. Math.,
{\textbf 53} (2000), 1--75.

\bibitem{Uh89}
Uhlenbeck, K., \emph{{H}armonic maps into Lie groups (classical
solutions of the chiral model)}, J. Diff. Geom., {\textbf 30}
(1989), 1--50.

\bibitem{Uhl92}
Uhlenbeck, K., \emph{{O}n the connection between harmonic maps and the self-dual Yang-Mills and the sine-Gordon equations}, J. Geom. Phys. \textbf{8} (1992), 283-316

\bibitem{Vi90}
Villarroel, J., \emph{{T}he inverse problem for Ward's system},
Stud. Appl. Math., {\textbf 83} (1990), 211--222.

\bibitem{Wa88}
Ward, R.S., \emph{{S}oliton solutions in an integrable chiral
model in $2+1$ dimensions}, J. Math. Phys., {\textbf 29} (1988),
386--389.

\bibitem{Wa90}
Ward, R.S., \emph{{C}lassical solutions of the chiral model,
unitons, and holomorphic vector bundles}, Comm. Math. Phys.,
{\textbf 128} (1990), 319--332.

\bibitem{Wa95}
Ward, R.S., \emph{{N}ontrivial scattering of localized solutions
in a $(2+1)$-dimensional integrable systems}, Phys. Letter A,
{\textbf 208} (1995), 203--208.

\bibitem{ZakMan73}
Zakharov, V.E., Manakov, S.V., \emph{{T}he theory of resonant interaction of wave
packets in non-linear media}, Sov. Phys. JETP \textbf{42} (1974), 842-850

\bibitem{ZaMi78}
Zakharov, V. E. and Mikhailov, A. V., \emph{{R}elativistically
invariant two dimensional models of fields theory which are
integrable by means of the inverse scattering problem method},
Sov. Phys. JETP, {\textbf 47} (1978), no. 6, 1017--1027.

\bibitem{ZakSha72}
Zakharov, V.E., Shabat, A.B., \emph{{E}xact theory of two-dimensional
self-focusing and one-dimensional of waves in nonlinear
media}, Sov. Phys. JETP \textbf{34}  (1972), 62-69

\bibitem{ZakSha79}
Zakharov, V.E., Shabat, A.B., \emph{{I}ntegration of non-linear equations of
mathematical physics by the inverse scattering method, II}, Funct. Anal.
Appl., \textbf{13} (1979), 166--174


\end{thebibliography}
\end{document}